\documentstyle[12pt]{amsart}
\newcommand{\ie}{{\it i.e.\/}\ }

\def\C{{\Bbb C}}
\def\Q{{\Bbb Q}}

\def\E{{\cal E}}
\def\O{{\cal O}}
\def\a{\alpha }

\def\si{\sigma }

\def\Ga{\Gamma }
\def\na{\nabla }
\def\om{\omega }
\def\r{\rho }

\def\Om{\Omega }
\def\th{\theta }
\def\p{\phi }

\long\def\comment#1\endcomment{}

\begin{document}

\title[Principal bundles admitting a holomorphic connection]{Principal
bundles admitting a holomorphic connection}
\author{Indranil Biswas}

\address{School of Mathematics, Tata Institute of Fundamental
Research, Homi Bhabha Road, Bombay 400005, INDIA}
\curraddr{Universit\'e de Grenoble 1, Institut Fourier,
100 rue des Maths, B.P. 74, 38402 Saint-Martin-d'H\`eres, FRANCE}
\email{indranil@@math.tifr.res.in and biswas@@puccini.ujf-grenoble.fr}
\thanks{Supported by the Acad\'emie des Sciences}
\date{}

\maketitle

\section{Introduction}

Let $E$ be a holomorphic vector bundle over a
compact K\"ahler manifold $M$, such that $E$ admits a holomorphic
connection compatible with its
holomorphic structure. The following question is perhaps
due to S. Murakami \cite{M1}, \cite{M2}, \cite{M3} : does there
exists a flat connection on $E$ compatible with the holomorphic
structure ?

The same question can be posed for a more general principal
$G$ bundle where $G$ is a complex Lie group. Murakami constructed
examples of torus bundles over torus which admit holomorphic
connection but do not admit compatible flat connection; thus
providing a negative answer to the general question.

Let $M$ be a compact connected K\"ahler manifold of complex
dimension $d$, and let $\om$
be a K\"ahler form on $M$. The degree
of a torsion-free ${\O}_M$-coherent sheaf $F$ is defined as
[Ko, Ch. V, (7.1)] $$\deg F \, := \, \int_Mc_1(F)\wedge
{\om}^{d-1} \leqno{(1.1)}$$
Consider the Harder-Narasimhan filtration
of the holomorphic tangent bundle of $M$ [Ko, page 174, Ch. V,
Theorem 7.15] : $$0 \, = \, {\E}_0 \, \subset \,
{\E}_1 \, \subset \, {\E}_2
\, \subset \, {\E}_3 \, \subset \, \ldots \, \subset \,
{\E}_{l-1} \, \subset \, {\E}_l \, = \, T \leqno{(1.2)}$$

Let $G$ be a connected affine algebraic reductive group over $\C$.
Let $P$ be a holomorphic principal $G$ bundle on $M$ (\ie the
transition functions of $P$ are holomorphic). Assume that $P$ admits
a holomorphic connection $D$ compatible with the holomorphic
structure. This means the following :
$D$ is a holomorphic $1$-form on $P$ with values in the Lie
algebra, $\frak g$, of $G$, such that $D$ is invariant for the
action of $G$ on $P$, and, when restricted to the fibers of
$P$, this form coincides
with the holomorphic Maurer-Cartan form. Using the natural
identification of the holomorphic tangent space of $P$ with
its real tangent space, the holomorphic connection $D$ gives
a $G$ connection on $P$.

A flat $G$ connection on $P$ is called compatible with
respect to the holomorphic structure if the $(1,0)$ component
of the connection form (which is a $\frak g$ valued $1$-form
on $P$) is actually a holomorphic $1$-form. This is equivalent
to the following : $M$ can be covered by open set $\{U_i\}$
such that over each $U_i$ the principal bundle $P$ admits
a holomorphic trivialization which is also constant (with
respect to the connection). Clearly the $(1,0)$ component
of such a flat connection compatible with the holomorphic
structure gives a holomorphic connection compatible with
the holomorphic structure.

Our aim here is to to prove the following theorem [Theorem 3.1] :

\medskip
\noindent {\bf Theorem A.}\,\ {\it If the K\"ahler manifold
$M$ satisfies the condition that $\deg (T/{\E}_{l-1}) \geq 0$, then
a holomorphic principal $G$-bundle $P$ on $M$ admitting
a compatible holomorphic connection is semistable. Moreover,
if $\deg (T/{\E}_{l-1}) >0$, then such a bundle $P$ actually
admits a compatible flat $G$-connection.}
\medskip

So, in particular, if $T$ is semistable with $\deg T >0$,
then a $G$-bundle $P$ on $M$ with a holomorphic connection
admits a compatible flat connection. In the case where
$G = GL(n, \C)$, the above result was proved in \cite{Bi}.

The author is grateful to D. Akhiezer for some very useful
discussions.

\section{Preliminaries}

We continue with the notation of Section 1.

Let $D$ be a holomorphic connection on the holomorphic principal
$G$ bundle $P$ compatible with its holomorphic structure.

Let $Ad(P)$ denote the vector bundle on $M$ associated
to $P$ for the adjoint action of $G$ on its Lie algebra
$\frak g$. The holomorphic structure of $P$ will induce
a holomorphic structure on the vector bundle $Ad(P)$.
Let ${\overline{\partial}}_P$ denote the differential
operator of order one defining the holomorphic structure of
$Ad(P)$.

Let ${\Om}^i_M$ denote the holomorphic vector bundle on $M$
given by the holomorphic $i$-forms.

The holomorphic connection $D$ on $P$ induces a holomorphic
connection on $Ad(P)$, which is again denoted by $D$. In other
words,
$$D \, : \, Ad (P) \, \longrightarrow \, Ad (P) \otimes {\Om}^1_M$$
is a first order operator satisfying
the Leibniz
condition $D(f.s) = \partial f . s + fD(s)$, where $f$ is a smooth
function, such that the curvature,
$(D+ {\overline{\partial}}_P)^2$, of
the connection $D+{\overline{\partial}}_P$ is a
holomorphic section of ${\Om}^2_M \otimes Ad (P)$. This condition
implies that the operator $D$ maps a holomorphic sections of
$Ad (P)$ to a holomorphic section of $Ad (P)\otimes {\Om}^1_M$.

A vector bundle $E$ on $M$ is called {\it stable} (resp.
{\it semistable}) if for any ${\O}_M$-coherent proper subsheaf
$0 \neq F \subset E$, with $E/F$ torsion-free, the following
condition holds [Ko, Ch. V, \S 7] :
$$ \mu (F) := \deg F/{\rm rank}F \, < \mu (E) :=
\deg E/{\rm rank} E \, ~ \, ~\, ({\rm resp.}\,\, \mu (F) \, \leq \,
\mu (E))$$ A vector bundle is called {\it quasistable} if it is
a direct sum of stable bundles of same $\mu$ (slope).

\medskip
\noindent {\bf Lemma 2.1.}\,\, {\it If $\deg (T/{\E}_{l-1}) \geq 0$
then the vector bundle $Ad (P)$ on $M$ is semistable. Moreover,
if $\deg (T/{\E}_{l-1}) > 0$, then $Ad (P)$ is quasistable.}
\medskip

The proof of this lemma is actually contained in \cite{Bi}.
However, to be somewhat self-contained, we will give some
details of the proof.

\medskip
\noindent {\underline {Proof of Lemma 2.1.}}\, We will first
prove that if $\deg (T/{\E}_{l-1}) \geq 0$ then the bundle
$Ad (P)$ is semistable.

Suppose $Ad (P)$ is not semistable. In that case it has a nontrivial
Harder-Narasimhan filtration. Let
$$0 \, = \, V_0\, \subset \, V_1 \, \subset \, V_2 \, \subset \,
\ldots \, \subset \, V_{n-1} \, \subset \,
V_n \, =\, Ad (P) \leqno{(2.2)}$$
be the Harder-Narasimhan filtration of the vector bundle $Ad (P)$.

We restrict the domain of the operator $D$ (which gives the
holomorphic connection on $Ad (P)$) to
the subsheaf $V_1$, and consider the induced operator
$$D_1 \, : \, V_1 \, \longrightarrow \,
(Ad(P)/V_1)\otimes {\Om}^1_M \leqno{(2.3)}$$ The Leibniz identity
for $D$ implies that $D_1(f.s) = f.D_1(s)$, \ie $D_1$ is
${\O}_M$-linear.

The largest semistable subsheaf (\ie the first nonzero term in the
Harder-Narasimhan filtration) of ${\Om}^1_M$ is the kernel
of the surjective homomorphism
$$q_{l-1} \, :\, {\Om}^1_M \, \longrightarrow \, {\E}^*_{l-1}$$
obtained by taking the dual of the inclusion in (1.2). Since
the tensor product of two semistable sheaves is again semistable
[MR, Remark 6.6 iii], the largest
semistable subsheaf of $(Ad (P)/V_1) \otimes {\Om}^1_M$ is
$$W \, := \, (V_2/V_1)\otimes {\rm kernel}(q_{l-1})\leqno{(2.4)}$$
The general formula for the degree of a tensor product gives
$\mu (W) = \mu (V_2/V_1) + \mu ({\rm kernel}(q_{l-1}))$.
{}From the property of the Harder-Narasimhan filtration we have that
$\mu (V_2/V_1) < \mu (V_1)$. The assumption that
$\deg (T/{\E}_{l-1}) \geq 0$ implies that
$$\mu ({\rm kernel}(q_{l-1})) \, = \, - \mu (T/{\E}_{l-1}) \, \leq \,0$$
{}From these we conclude that
$$\mu (V_1) \, > \, \mu (W) \leqno{(2.5)}$$
Let $E$ denote the image of $D_1$ defined in (2.3). Assume
that $E$ is not the zero sheaf. Since
$V_1$ is semistable and $E$ is a quotient of $V$, we have
$\mu (E) \geq \mu (V_1)$. On the other hand, since $W$ is the
largest semistable subsheaf of $(Ad (P)/V_1)\otimes {\Om}^1_M$,
we have $\mu (E) \leq \mu (W)$. But this contradicts (2.5). So
$E$ must be the zero sheaf, \ie  the homomorphism
$D_1$ in (2.3) must be the zero homomorphism. Thus we obtain that the
subsheaf $V_1$ is invariant under the connection $D$ on $Ad (P)$.

Any ${\O}_M$ coherent sheaf
with a holomorphic connection is a locally free ${\O}_M$ module
[B, p. 211, Proposition 1.7]. (This proposition in \cite{B}
is stated for integrable connections ($D$-modules), but the proof
uses only the Leibniz rule which is valid for a holomorphic
connection.)

Using the Chern-Weil construction of characteristic classes it is
easy to show that for a vector bundle $E$ on $M$, equipped with
a holomorphic connection, any Chern class $c_i(E) \in H^{2i}(M,\Q)$,
$i\geq 1$, vanishes.

So we have $\deg V_1 = 0 = \deg Ad (P)$. If $V_1 \neq Ad(P)$
then $\mu (V_1) > \mu (Ad (P))$. So $Ad (P) = V_1$, \ie
the bundle $Ad (P)$ is semistable.

If $Ad (P)$ is stable then obviously
it is quasistable. Suppose that $Ad(P)$ is not stable.
Then there is a filtration [K, Ch. V, \S 7,
Theorem 7.18] $$0 \,= \,W_0 \, \subset \, W_1 \, \subset \,W_2 \,
\subset
\, \ldots \, \subset \, W_{n-1} \, \subset \, W_n \,
= \, Ad (P) \leqno{(2.6)}$$
such that $W_i/W_{i-1}$, $1\leq i \leq n$, is a stable
sheaf with $\mu (W_i/W_{i-1}) = \mu (Ad (P))$.

Now, as done in the proof of Proposition 3.4 of \cite{Bi},
using the given condition that $\deg (T/{\E}_{l-1})$ is
strictly positive, it is
possible to show that the filtration (2.6) splits, \ie
$$W_i \, = \, W_{i-1}\oplus (W_i/W_{i-1}) \leqno{(2.7)}$$
where $W_i/W_{i-1}$ is a locally free sheaf. The equality
(2.7) proves that $$Ad (P)\, =\,  \sum_{i=1}^n W_i/W_{i-1}$$
This completes the proof of the lemma. $\hfill{\Box}$
\medskip

In the next section we will use Lemma 2.1 to construct
a flat connection on $P$.

\section{Existence of a flat connection}

Let $G$ be a connected affine algebraic reductive group over $\C$.
Let $P$ be a holomorphic principal $G$ bundle over $M$.

Let $i : U \longrightarrow M$ be the inclusion of an open subset
such that the complement $X-i(U)$ is an analytic subset of $X$ of
codimension at least two. For a ${\O}_U$ coherent sheaf $F$
on $U$, the direct image $i_*F$ is a ${\O}_M$ coherent sheaf.
The degree of $F$ is defined to be the degree of $i_*F$.

We will recall the definition of (semi)stability of $P$
([RR, Definition 4.7]).
Let $U\subset X$ be an open subset
with $X-U$ being an analytic set of codimension at least two.
Let $Q$ be a parabolic subgroup of $G$, and let $P'$
be a reduction of the structure group to $Q$
of the restriction of the principal $P$ to the open set $U$. The
principal bundle $P$ is said to be {\it stable}
(resp. {\it semistable}) if for any such $P'$, the degree of a
line bundle on $U$ associated to $P$ for any character
$\chi$ on $Q$ dominant with respect to a Borel subgroup
contained in $Q$, is strictly negative (resp. nonpositive).

Let $D$ be a holomorphic connection on $P$ (defined in Section 1).
The following theorem is a generalization of Lemma 2.1.

\medskip
\noindent {\bf Theorem 3.1.}\, {\it If $\deg (T/{\E}_{l-1}) \geq 0$
then $P$ is semistable. Moreover, if $\deg (T/{\E}_{l-1}) >0$,
then $P$ admits a flat $G$-connection compatible with
the holomorphic structure.}
\medskip

\noindent {\bf Proof.}\, Let $P' \subset P$ be a reduction
of structure group of $P$ to a maximal parabolic subgroup
$Q\subset G$. This reduction is given by a section, $\si$,
of the fiber bundle
$${\r} \, : \, P/Q \, \longrightarrow \, U$$ Let $T_{\rm
rel}$ denote the relative tangent bundle for the map $\r$.

{}From Lemma 2.1 of \cite{R} it
follows that in order to check that $P$ is semistable it is enough
to show that $\deg ({\si}^*T_{\rm rel}) \geq 0$.

The reduction $P'\subset P$ gives an injective homomorphism
$Ad(P') \longrightarrow Ad(P)$ of adjoint bundles on $M$.
The bundle ${\si}^*T_{\rm rel}$ on $M$ is the quotient
bundle $Ad(P)/Ad(P')$.

{}From Lemma 2.1 we know that if $\deg (T/{\E}_{l-1}) \geq 0$, the
adjoint bundle $Ad(P)$ is semistable. Since $G$ is reductive,
the Lie algebra $\frak g$ admits a nondegenerate
$G$ invariant bilinear form. This implies that $Ad(P) = Ad(P)^*$.
Hence $\deg Ad(P) =0$. Now the semistabilty of $Ad(P)$
implies that $\deg (Ad(P)/Ad(P')) \geq 0$. This proves that
the principal bundle $P$ is semistable.

If $\deg (T/{\E}_{l-1}) > 0$ then from Lemma 2.1 we know that
$Ad (P)$ is quasistable. So from the main theorem of
\cite{UY} it follows that the vector bundle $Ad(P)$ admits
a Hermitian-Yang-Mills connection. We will denote this
Hermitian-Yang-Mills connection
by $\na$. This connection $\na$ is unique (though the
Hermitian-Yang-Mills metric is not unique), and it is irreducible
if and only if $Ad(P)$ is stable.

We want to show that this connection $\na$ induces a connection
on the principal $G$ bundle $P$.

Let $Z_0$ denote the connected component of the
center of $G$ containing the identity element
(the center has finitely many components).
Define $$G_0 \, = \, G/Z_0$$ which is a semisimple group.
The group $G_0$ acts on the Lie algebra $\frak g$ (of $G$)
by conjugation and gives an homomorphism
$$\th \, : \, G_0 \, \longrightarrow \,GL({\frak g})
\leqno{(3.2)}$$ which has a finite group as the kernel.

{}From a theorem of Chevalley [H, Theorem 11.2] we know that
there is a linear representation of the group $GL ({\frak g})$
$$\p \, : \, GL({\frak g}) \, \longrightarrow \, GL(V)
\leqno{(3.3)}$$
in a vector space $V$ over $\C$ and a line $L$ in $V$ such that
$$\p \circ \th (G_0) \, =
\, \{g \in GL({\frak g}) \, \vert \,~\, {\p}(g)(L) = L \}$$
Since $G_0$ is semisimple, it does not have
any nontrivial character. This implies that $\p\circ \th (G_0)$
fixes the line $L$ point-wise. Let $0\neq v \in L$ be a nonzero
vector. So the isotropy subgroup of $v$ for the action
of $GL({\frak g})$ on $V$ is precisely ${\p \circ\th}(G_0)$. It
is not difficult to see that the homomorphism $\p$ can be
so chosen that it maps the center of the Lie algebra of
$GL({\frak g})$ into the center of the Lie algebra of $GL(V)$.
We will choose $\p$ such that it satisfies this condition.

Let $q : G \longrightarrow G_0$ denote the obvious projection.
Let $P(G_0)$ denote the principal $G_0$ bundle on $M$ obtained by
extending the structure group of $P$ to $G_0$ using the
homomorphism $q$. The vector bundle $Ad(P)$, which can be
identified with a principal $GL({\frak g})$ bundle, is obtained by
extending the structure group of $P(G_0)$ to $GL({\frak g})$
using the homomorphism $\th$ defined in (3.2).

Using the homomorphism $\p$ we may extend the structure
group of $Ad(P)$ to $GL(V)$. The vector bundle on $M$ associated
to this principal $GL(V)$ bundle for the natural action of
$GL(V)$ on the vector space $V$ will be denoted by $E$.

The Hermitian-Yang-Mills metric on the vector bundle $Ad(P)$
gives a reduction of the structure group of $Ad(P)$ to
$U({\frak g})$, a maximal compact subgroup of $GL ({\frak g})$.
Since the image $\p (U({\frak g}))$ is a compact
subgroup of $GL(V)$, it is contained in some
maximal compact subgroup of $GL(V)$. So the reduction of $Ad(P)$
to $U({\frak g})$ will induce a reduction of $E$ to a maximal
compact subgroup of $GL(V)$ (\ie the vector bundle
$E$ will be equipped with a hermitian metric) such that the
connection on $E$ obtained by extending the connection
$\na$ (on $P({\frak g})$) is the hermitian connection (for the
hermitian metric on $E$).
Since the metric on $Ad(P)$ is a
Hermitian-Yang-Mills metric, the metric on $E$ obtained above
is also a Hermitian-Yang-Mills metric. Indeed, the
Hermitian-Yang-Mills condition of the connection on $Ad(P)$
implies that the curvature is a $2$-form on $M$ with values in the
center of the endomorphism bundle $Ad(P)^*\otimes Ad(P)$. Since $\p$
maps the center of the Lie algebra of $GL({\frak g})$ into the the
center of the Lie algebra of $GL(V)$, the induced connection on
$E$ is a Hermitian-Yang-Mills connection. Let ${\na}'$ denote the
Hermitian-Yang-Mills connection on $E$ obtained this way.

Since ${\th}(G_0)$ fixes the vector $v$, this vector $v$ will give
a nowhere zero section of $E$ (since $E$ is obtained by extending
the structure group of $P$ to $GL(V)$); let $s$
denote this section of $E$.

The holomorphic connection $D$ on $P$ will induce a holomorphic
connection on $E$. As we noted in Section 2, using
the Chern Weil construction it is easy to see that the existence
of a holomorphic connection on $E$ implies that any Chern class,
$c_i(E)$, $i\geq 1$, vanishes.

{}From [Ko, Ch. IV, Corollary 4.13] it follows that
Hermitian-Yang-Mills connection ${\na}'$
is a flat connection. Since ${\na}'$ is a flat unitary
connection, any holomorphic section
of $E$ must be a flat section for the connection ${\na}'$. Indeed,
the Laplacian of a flat unitary connection operator is twice the
Laplacian of the Dolbeault operator. In particular, the space of
harmonic sections for these two Laplacians coincide. Since any
holomorphic section is a harmonic section for the Dolbeault
Laplacian, it must be a flat section. So, in particular,
the section $s$ is flat.

Let $P({\frak g})$ denote the principal $GL({\frak g})$
obtained by extending the structure group of $P$ to
$GL({\frak g})$ using the homomorphism $\th\circ q$. Since
$P({\frak g})$ is an extension, the fiber bundle,
$P({\frak g})/G$, with fiber $GL({\frak g})/G$ has a natural
section (which gives the reduction of the structure
group of $P({\frak g})$ to $G$). Let $\a$ denote this section
of $P({\frak g})/G$.

Since $\th (G_0)$ fixes $v$ for the homomorphism $\p$ in (3.3),
we have an embedding of the fiber bundle $P({\frak g})/G$
in the total space of $E$ (given by the orbit of $v$ for
the action of $GL({\frak g})$ using $\p$).

It is easy to see that the image of the section $\a$ by the
above embedding of $P({\frak g})/G$ in $E$
is precisely the section $s$.

Now, since $s$ is a flat section for the connection ${\na}'$,
the connection $\na$ induces a $G_0$ connection on the principal
$G_0$ bundle $P(G_0)$ on $M$ as follows : Let $$p \, :\,
P(G_0) \, \longrightarrow \, P({\frak g})$$ denote the
holomorphic map induced by $\th$ in (3.2). Take $x \in p (P(G_0))$,
and let $v \in T_xP({\frak g})$ be a horizontal vector for the
connection $\na$ on $P({\frak g})$. Let $w$ be the image
of $v$ by the differential of the map
$$P({\frak g})\, \longrightarrow \, E$$
induced by the homomorphism $\p$ in (3.3). Since $s$ is flat,
the tangent vector $w$ lies in the submanifold of the total
space of $E$ given by the section $s$. But this implies that
the tangent vector $v$ lies in the image of $TP(G_0)$ under
the map given by the differential of $p$. Thus the connection
$\na$ on $P({\frak g})$ induces a connection on the principal
bundle $P(G_0)$ with structure group $G_0$. We will call
this connection on $P(G_0)$ as ${\na}_0$. Since $\na$ is a
flat connection, ${\na}_0$ is also flat.

The commutator subgroup
$G' \, := \,[G,G] \, \subset \,G$ is a semisimple group, and
the restriction of $q$ to $G'$ is a surjective homomorphism with a
finite kernel. So their Lie algebras are isomorphic. So $$G \,=
Z_0.G' \leqno{(3.4)}$$
with a finite intersection $\Ga := Z_0\cap G'$.

The abelian Lie group $Z_0/\Ga$ is a product of copies of
${\C}^*$, since $G$ is assumed to be affine. Let
$$f \, :\, G \,\longrightarrow \, Z_0/\Ga$$ denote the
obvious projection (obtained from (3.4)).

Let $P(f)$ denote the principal $Z_0/\Ga$ bundle on $M$ obtained
by extending the structure group of $P$ using the homomorphism
$f$.

The holomorphic connection $D$ on $P$ induces a holomorphic
connection on $P(f)$, which we will denote by $D(f)$.

Any holomorphic line bundle on $M$ admitting a holomorphic
connection actually admits a compatible flat connection. Indeed,
if $\partial$ is a holomorphic connection on a holomorphic
line bundle $L$ whose holomorphic structure is given by the
operator $\overline\partial$, then the curvature
$(\partial +{\overline{\partial}})^2$ is a holomorphic $2$-form
which is exact (since the cohomology class represented
by it is of the type $(1,1)$).  So it is of the form
$\partial \beta$, where $\beta$ is a $(1,0)$-form. The new
connection $$\partial -\beta + {\overline{\partial}}$$ on $L$ is a
flat connection compatible with the holomorphic structure.

Recall that $Z_0/\Ga$ is a product of copies of ${\C}^*$. In
view of the above remark, the the existence of the holomorphic
connection $D(f)$ implies that the principal $Z_0/\Ga$ bundle
$P(f)$ admits a flat connection. Let ${\na}_1$ be a flat
connection on $P(f)$.

Since the exact sequence of the Lie algebras
$$0 \, \longrightarrow \, {\frak z}_0 \,\longrightarrow \,
{\frak g} \,\longrightarrow \,{\frak g}_0 \,\longrightarrow\,
0$$ has a natural splitting (given by the Lie algebra of $G'$),
the two flat connections ${\na}_0$ and ${\na}_1$ combine
together to induce a flat $G$ connection on $P$ as follows :
The horizontal subspace of the tangent space at a point $p \in P$
is defined to be the intersection of the inverse images of the
horizontal subspaces of $P(G_0)$ and $P(f)$ (horizontal
subspaces for the
flat connections ${\na}_0$ and ${\na}_1$ respectively) for the
obvious projections of $P$ onto $P(G_0)$ and $P(f)$ respectively.
The integrability of ${\na}_0$ and ${\na}_1$ will imply that
the connection on $P$ obtained above is actually flat.
This completes the proof of the theorem. $\hfill{\Box}$



\begin{thebibliography}{99}

\bibitem[B]{B} A. Borel et al. : Algebraic $D$-modules.
Perspectives in Mathematics, Vol. 2 (Ed. J. Coates, S. Helgason),
Academic Press, 1987.

\bibitem[Bi]{Bi} I. Biswas : On Harder-Narasimhan filtration
of the tangent bundle. To appear in Comm. Anal. Geom.

\bibitem[H]{H} J. Humphreys : Linear algebraic groups. Graduate
Texts in Math. 21, Springer-Verlag, New York Heidelberg
Berlin, 1975.

\bibitem[Ko]{Ko} S. Kobayashi : Differential geometry of complex
vector bundles. Publications of Math. Soc. of Japan, Iwanami
Schoten Pub. and Princeton University Press, 1987.

\bibitem[M1]{M1} S. Murakami : Sur certains espaces fibr\'es
principaux admettant des connexions holomorphes. Osaka Math. Jour.
{\bf 11} (1959) 43-62.

\bibitem[M2]{M2} S. Murakami : Sur certains espaces fibr\'es
principaux holomorphes dont le groupe est ab\'elien connexe.
Osaka Math. Jour. {\bf 13} (1961) 143-167.

\bibitem[M3]{M3} S. Murakami : Harmonic connections and their
applications. SEA Bull. Math. Special Issue (1993) 101-103.
World Sci. Pub. Company.

\bibitem[MR]{MR} V. Mehta, A. Ramanathan : Semistable sheaves
on projective varieties and their restriction to curves. Math.
Ann. {\bf 258} (1982) 213-224.

\bibitem[RR]{RR} S. Ramanan, A. Ramanathan : Some remarks on
the instability flag. T\^ohoku Math. J. {\bf 36} (1984) 269-291.

\bibitem[R]{R} A. Ramanathan : Stable principal bundles on a
compact Riemann surface. Math. Ann. {\bf 213} (1975) 129-152.

\bibitem[UY]{UY} K. K. Uhlenbeck, S. T. Yau : On the existence
of Hermitian-Yang-Mills connections in stable vector bundles.
Comm. Pure Appl. Math. {\bf 39} (1986) 257-293.

\end{thebibliography}
\end{document}